\documentclass[11pt]{article}
\usepackage[utf8]{inputenc}
\usepackage[colorlinks=true,urlcolor=blue,
citecolor=red,linkcolor=blue,linktocpage,pdfpagelabels,bookmarksnumbered,bookmarksopen]{hyperref}
\usepackage[numbers]{natbib}

\usepackage{amsmath}
\usepackage{amssymb}
\usepackage{latexsym}
\usepackage{enumitem}
\usepackage{mathrsfs}
\usepackage{comment}
\usepackage{color}
\usepackage{bbm}

\newtheorem{assumption}{Assumption}

\def\qed{ \ \vrule width.2cm height.2cm depth0cm\smallskip}
\newenvironment{proof}{\noindent {\bf Proof.\/}}{$\qed$\vskip 0.1in}

\newcommand{\ol}{\overline}

\newcommand{\ba}{\begin{array}}
\newcommand{\ea}{\end{array}}
\newcommand{\be}{\begin{equation}}
\newcommand{\ee}{\end{equation}}
\newcommand{\bea}{\begin{eqnarray}}
\newcommand{\eea}{\end{eqnarray}}
\newcommand{\beaa}{\begin{eqnarray*}}
\newcommand{\eeaa}{\end{eqnarray*}}

\def\dbE{\mathbb{E}}
\def\dbF{\mathbb{F}}

\def\dbH{\mathbb{H}}

\def\dbL{\mathbb{L}}

\def\dbP{\mathbb{P}}
\def\dbR{\mathbb{R}}

\def\dbV{\mathbb{V}}
\def\dbX{\mathbb{X}}
\def\dbY{\mathbb{Y}}

\def\a{\alpha}
\def\b{\beta}

\def\d{\delta}
\def\e{\varepsilon}

\def\t{\tau}

\def\th{\theta}

\def\o{\omega}

\def\D{\Delta}
\def\Th{\Theta}

\def\O{\Omega}
\def\cA{{\cal A}}

\def\cE{{\cal E}}
\def\cF{{\cal F}}

\def\cH{{\cal H}}
\def\cI{{\cal I}}

\def\no{\noindent}

\def\ss{\smallskip}
\def\ms{\medskip}
\def\bs{\bigskip}
\def\q{\quad}
\def\qq{\qquad}

\def\pa{\partial}
\def\cd{\cdot}
\def\cds{\cdots}

\def\bx{{\bf x}}

\def\tr{\hbox{\rm tr}}

\def\qed{ \hfill \vrule width.25cm height.25cm depth0cm\smallskip}

\newcommand{\basa}{\begin{assumption}}
\newcommand{\easa}{\end{assumption}}

\newcommand{\bas}{\begin{assum}}
\newcommand{\eas}{\end{assum}}

\def\limsup{\mathop{\overline{\rm lim}}}

\def\pa{\partial}

 \def\cd{\cdot}
\def\cds{\cdots}

\def\tr{\hbox{\rm tr$\,$}}

\def\dis{\displaystyle}

\def\bx{{\bf x}}

\def\1{{\bf 1}}

\def\:{\!:\!}
\def\reff#1{{\rm(\ref{#1})}}
\def \proof{{\noindent \bf Proof\quad}}

at 9pt

\newtheorem{thm}{Theorem}[section]
\newtheorem{lem}[thm]{Lemma}

\newtheorem{prop}[thm]{Proposition}
\newtheorem{rem}[thm]{Remark}
\newtheorem{eg}[thm]{Example}
\newtheorem{defn}[thm]{Definition}
\newtheorem{assum}[thm]{Assumption}

\numberwithin{equation}{section}

\begin{document}

\title{\bf  Set Values of Dynamic Nonzero Sum Games and Set Valued Hamiltonians}
\author{, Jianfeng Zhang}

\author{
Bixing Qiao\thanks{\noindent Department of
Mathematics, University of Southern California, Los Angeles, 90089; email: bqiao@usc.edu.}
 ~ and ~ Jianfeng Zhang \thanks{ \noindent Department of
Mathematics, University of Southern California, Los Angeles, 90089;
email: jianfenz@usc.edu. This author is supported in part by NSF grant  \#DMS-2205972. }}

\date{\today}
\maketitle

\begin{abstract}
It is well known that the (unique) value of a stochastic control problem or a two person zero sum game under Isaacs condition can be characterized through a PDE driven by the Hamiltonian. Our goal of this paper is to extend this classical result to nonzero sum games, which typically have multiple Nash equilibria and multiple values. Our object is the set value of the game, which roughly speaking is the set of values over all equilibria and thus is by nature unique. We shall introduce set valued Hamiltonians and characterize the set value of the game through backward SDEs driven by appropriate selectors of the set valued Hamiltonians, where the selectors are typically path dependent. When the set valued Hamiltonian is a singleton, our result covers the standard control problem and two person zero sum game problem under Isaacs condition.
 \end{abstract}

 \bs
\no{\bf Keywords.} Nonzero sum games, Nash equilibria, Set values, Hamiltonian

\bs

\no{\it 2020 AMS Mathematics subject classification:} 91A15, 49L12, 60H10

\vfill\eject
\section{Introduction}
\label{sect-introduction}

The general theory of games can be traced back to von Neumann-Morgenstern \cite{VNM} and Nash \cite{Nash2}, and we refer to the textbook Fudenberg-Tirole \cite{FudenbergTirole} for an excellent exposition. In this paper, we focus on the basic notion Nash equilibrium, which will often be simply referred to as equilibrium in the paper, for nonzero sum stochastic differential games. Unlike stochastic control problems where the unique value is always well defined, a game typically has multiple values induced by multiple equilibria. This feature hampers the efforts to study or even just to define the value of general games, especially in a dynamic setting.  Indeed, there have been very few works on dynamical approaches for general nonzero sum games. The works Feinstein-Rudloff-Zhang \cite{feinstein2020dynamic} and Iseri-Zhang \cite{iseri2024set} proposed to study the set value of games, which roughly speaking is the set of values over all equilibria. See also Abreu-Pearce-Stacchetti \cite{APS},  Cardaliaguet-Plaskacz \cite{Cardaliaguet2003}, Cardaliaguet-Quincampoix-Saint Pierre \cite{CQS},  Feinstein \cite{Feinstein2022}, Lazrak-Zhang \cite{ LazrakZhang}, Sannikov \cite{Sannikov}, and Zhang \cite{ZhangEfficiency} for some related works. The set value is by definition unique, and it shares many nice properties of the standard value function of stochastic control problems, most notably the dynamic programming principle, also called time consistency.

The dynamic programming can be viewed as a bridge between local (in time) and global equilibrium property, but in the integration or expectation form.  For standard stochastic control problems, an alternative bridge is through the Hamiltonian in the differentiation form. Indeed, based on the dynamic programming principle, the dynamic value function of the control problem is typically characterized as the unique viscosity solution of the Hamilton-Jacobi-Bellman equation driven by the Hamiltonian, c.f. Fleming-Soner \cite{FlemingSoner} and Yong-Zhou \cite{yong_stochastic_1999}. Similar results hold for two person zero sum games under Isaacs condition, c.f. Fleming-Souganidis \cite{Fleming1989} and Hamadene-Lepeltier \cite{HAMADENE1995259}. The goal of this paper is to extend these results to nonzero sum games with multiple equilibria and more importantly with multiple values. To be specific, we shall introduce the Hamiltonian for the nonzero sum game, which also takes set value and thus is called set valued Hamiltonian, and then we characterize the set value of the game thorough its set valued Hamiltonian. We shall remark that the set valued Hamiltonian has much simpler structure than the set value of the game and is a lot easier to analyze.

Our main results are in three folds. First, we consider the case that the set valued Hamiltonian is a singleton and thus its (unique) element is naturally called the (vector valued) Hamiltonian of the game. We show that the set value of the game is also a singleton and the unique game value is characterized by the system of HJB equations driven by the Hamiltonian. This result
covers the standard stochastic control problem, which can be viewed as a game with one player, and two person zero sum game problem under Isaacs condition. We remark that the Isaacs condition exactly means that set valued Hamiltonian of the two person zero sum game is a singleton, and when the Isaacs condition fails, the set valued Hamiltonian is empty and thus the game has no value.  We shall also emphasize that, as in \cite{feinstein2020dynamic, iseri2023set} as well as in Buckdahn-Cardaliaguet-Rainer \cite{Buckdahn2004} and Mertens-Sorin-Zamir \cite[Chapter VII.4]{MSZ}, both for the game and for the Hamiltonian our set value consists of all values which are limits of values at approximate equilibria, rather than the raw set value which consists of all values at true equilibria. In particular, the set value may not be empty even if there is no equilibrium, see e.g. Frei-dos Reis \cite{DF} for a counterexample. This is consistent with the standard literature on control problems and two person zero sum game problems where one studies their values without requiring the existence of optimal controls or saddle points.

We next study the case that the set valued Hamiltonian is Lipschitz continuous and is separable, in the sense that it contains a dense subset of uniformly Lipschitz continuous selectors. In this case, we may characterize the set value of the game through backward SDEs driven by  Lipschitz continuous selectors of the set valued Hamiltonian. This result effectively transforms the game problem into a control problem, with the selectors as the control, and it provides a powerful and convenient tool for further analysis of the game. We emphasize that, while the set valued Hamiltonian is state dependent, we have to allow the selectors to be path dependent. This is consistent with the observation in \cite{feinstein2020dynamic} that the dynamic programming of the game set value  requires to consider path dependent equilibria. We thus use backward SDE of Pardoux-Peng \cite{PardouxPeng} which allows for path dependence. Alternatively, we may also consider path dependent PDEs, see e.g. Zhang \cite[Chapter 11]{Zhangbook}.

We finally study the case that the set valued Hamiltonian is merely continuous. In this case, we are not able to express the game set value through true selectors of the set valued Hamiltonian. Instead, we characterize it as the set of limit values of BSDEs driven by approximate selectors, in the same spirit as we define set values. We shall remark though, game values are quite sensitive to small perturbations of game parameters, and thus the continuity of the set valued Hamiltonian is not a trivial requirement. It will be very interesting to relax this regularity requirement and to extend our results further.

At this point we should mention that there have been many works in the literature which used Hamiltonians, either the unique one or some appropriate selector, to construct equilibria of nonzero sum games,  see e.g. Bensoussan-Frehse  \cite{Bensoussan2000}, El Karoui-Hamadene \cite{ELKAROUI2003145}, Espinosa-Touzi \cite{EspinosaTouzi}, Friedman \cite{FRIEDMAN197279},     Hamadene-Lepeltier-Peng \cite{zbMATH01066318}, Hamadene-Mu  \cite{HAMADENE2014699, HamadeneMu2, HamadeneMu}, and Uchida \cite{Uchida}. We should note though that all these works talk about true equilibria, both for the Hamiltonian and for the game, rather than the limit of approximate equilibria values. More importantly, these works did not study the opposite direction that all values of the game can be constructed in this way. We instead characterize the whole set value through the set valued Hamiltonians. 

We would also like to mention another highly related direction on set valued PDEs. In the present paper, while the game values and the Hamiltonians are set valued, the involved systems of HJB equations and BSDEs are vector valued and thus are in the standard sense. In the contexts of multivariate control problems,  Iseri-Zhang \cite{iseri2023set} derived a set value HJB equation, where the set value itself becomes the solution of the PDE. The ongoing work Iseri-Zhang \cite{IseriZhangGamePDE} extends this result to set values of nonzero sum games.  Our set valued Hamiltonian plays a crucial role in the set valued PDE there. 

 The rest of the paper is organized as follows. In \S\ref{sect-setvalue} we introduce the problem. In \S\ref{sect-continuous} we study the general case where the set valued Hamiltonian is only continuous. In \S\ref{sect-Lipschitz} we study the case that the set valued Hamiltonian is Lipschitz continuous, which includes the singleton case. Finally in Appendix \ref{sect-appendix} we complete some postponed proofs.

\section{Set values and set valued Hamiltonians for games}
\label{sect-setvalue}
Let $(\O, \cF, \dbP)$ be a probability space, $B$ a standard $d$-dimensional Brownian motion, and $\dbF = \dbF^B$. We consider an $N$-player nonzero sum stochastic differential game on $[0, T]$ with drift controls. As standard in the literature, we consider closed loop controls and thus use weak formulation. That is, we shall fix the state process $X$ and the players control its law:
\beaa
X \equiv B.
\eeaa

For each $i=1,\cds, N$, let $\cA_i$ denote the set of admissible controls $\a^i$ for  player $i$ which are $\dbF$-progressively measurable and $A_i$-valued processes, where $A_i$ is a domain in certain Euclidean space. Denote $A := A_1\times \cds\times A_N$,  $\cA := \cA_1\times \cds\cA_N$, and $\a := (\a^1, \cds, \a^N)$. We emphasize again that, since we are using the canonical space of $X$, our controls are closed loop type rather than open loop type. Our game involves progressively measurable data:
\beaa
b, f: [0, T]\times \dbR^d\times A\to \dbR^N,\q g:  \dbR^d\to \dbR^N.
\eeaa
 Throughout the paper, the following assumption will always be in force.

\begin{assum}
\label{assum-standing}
$b, f, g$ are bounded and uniformly continuous in all variables.
\end{assum}
We remark that the boundedness assumption, especially that of $f$ and $g$, can be weakened. However, in that case the boundedness and  compactness of the set values may fail and the corresponding statements in the paper should be modified appropriately.

For any $\a\in \cA$, define
\bea
\label{PaJa}
\left.\ba{c}
\dis  M^\a_T := \exp\Big(\int_0^T b(t, X_t, \a_t) \cd dB_t - {1\over 2}\int_0^T |b(t, X_t, \a_t)|^2 dt\Big),\q \dbP^\a := M^\a_T d\dbP,\ss\\ 
\dis J_i(\a) := \dbE^{\dbP^\a} \Big[g_i(X_T) + \int_0^T f_i(t, X_t, \a_t)dt\Big],\q i=1,\cds, N,
\ea\right.
\eea
and $J(\a):= (J_1(\a),\cds, J_N(\a))$. We note that, in many applications $f_i$ may depend only on $\a^i$ instead of the whole vector $\a$.  

We next define Nash equilibrium and approximate Nash equilibrium of the game.

\begin{defn}
\label{defn-NE}
(i) We say $\a^*\in \cA$ is a Nash Equilibrium, denoted as $\a^*\in \cE$, if
\beaa
J_i(\a^*) \ge J_i(\a^{*,-i}, \a^i),\q\forall i=1,\cds, N, ~ \a^i \in \cA_i.
\eeaa

\no(ii) For any $\e>0$, we say $\a^\e\in \cA$ is an $\e$-Nash Equilibrium, denoted as $\a^\e\in \cE_\e$, if
\beaa
J_i(\a^\e) \ge J_i(\a^{\e,-i}, \a^i)-\e,\q\forall i=1,\cds, N, ~ \a^i \in \cA_i.
\eeaa
\end{defn}
Here $(\a^{*, -i},\a^i)$ is the vector by replacing $\a^{*, i}$ in $\a^*$ with $\a^i$, similarly for $(\a^{\e,-i}, \a^i)$. As in \cite{feinstein2020dynamic}, we define the raw set value $\dbV_0$ and the set value $\dbV$ of the game:
\bea
\label{setV}
\left.\ba{c}
\dis \dbV_0 := \big\{J(\a^*): \a^*\in \cE\big\} \subset \dbR^N;\ss\\
\dis \dbV := \bigcap_{\e>0} \dbV_\e \subset \dbR^N,\q \dbV_\e := \big\{y\in \dbR^N: |y- J(\a^\e)|< \e ~\mbox{for some}~\a^\e\in \cE_\e\big\}.
\ea\right.
\eea

\begin{rem}
\label{rem-raw}
(i) As explained in \cite{feinstein2020dynamic}, the set value $\dbV$, rather than the raw set value $\dbV_0$, is the natural counterpart of the value of standard stochastic control problems as well as the game value of two person zero sum game problems under the Isaacs condition. In particular, the definition of set value $\dbV$ does not require the existence of Nash equilibria. In fact, it is possible that $\dbV \neq \emptyset = \dbV_0$.

(ii) One can similarly define the dynamic set value function $\dbV(t,x)$ for the game on $[t, T]$ with initial condition $X_t=x$. One of the main results in \cite{feinstein2020dynamic} is that $\dbV(t,x)$ satisfies the dynamic programming principle in appropriate sense. The main goal of this paper is to characterize $\dbV=\dbV(0,0)$ through set valued Hamiltonian $\dbH$ which will be defined next. Similarly we can characterize $\dbV(t,x)$ through $\dbH$, and this characterization can serve as an alternative proof for the dynamic programming principle of $\dbV(t,x)$.     
\end{rem}

We now turn to the Hamiltonians. First, it is clear that 
\beaa
J_i(\a) = Y^{\a,i}_0,
\eeaa
 where $(Y^{\a,i}, Z^{\a,i})$ is the solution to the following linear BSDE taking values in $\dbR\times \dbR^d$:
\bea
\label{YZa}
\left.\ba{c}
\dis Y^{\a,i}_t = g_i(X_T) + \int_t^T h_i(s, X_s, Z^{\a,i}_s, \a_s) ds - \int_t^T Z^{\a,i}_s \cd dB_s,\ss\\
\dis \mbox{where}\q h_i(t,x,z^i, a) := f_i(t,x, a) + b(t,x,a) \cd z^i.
\ea\right.
\eea
Note that $h$ has the following regularity: for some modulus of continuity function $\rho_0$,
\bea
\label{hreg}
\big|h(\tilde t, \tilde x, \tilde z, \tilde a) - h(t,x,z,a)\big| \le C|\tilde z-z| + C[1+|z|] \rho_0\big(|\tilde t-t|+|\tilde x-x|+|\tilde a-a|\big).
\eea
Introduce the extended state space 
\bea
\label{dbX}
\dbX :=  [0, T]\times \dbR^d \times \dbR^{dN},\q\mbox{with its elements denoted as}~ \th = (t,x,z).
\eea
For such an $\th$, denote $\th^i := (t,x, z^i)\in [0, T]\times \dbR^d \times \dbR^{d}$, $i=1,\cds, N$, and set
 $h(\th, a):= \big(h_1(\th^1, a),\cds, h_N(\th^N, a)\big)$. Note that the mapping $a\in A\mapsto h(\th,a)\in \dbR^N$ is a static $N$-player game. Define Nash equilibrium and $\e$-Nash equilibrium of this static game in an obvious manner:
\bea
\label{hNE}
\left.\ba{c}
\dis a^* \in \cE(\th)\q \mbox{if}\q h_i(\th^i, a^*)  \ge  h_i(\th^i, (a^{*,-i}, \tilde a_i)),\q \forall i=1,\cds, N,~  \tilde a_i\in A_i;\ss\\
\dis a^\e \in \cE_\e(\th)\q\mbox{if}\q h_i(\th^i, a^\e) \ge h_i(\th^i, (a^{\e,-i}, \tilde a_i))-\e,\q \forall i=1,\cds, N,~  \tilde a_i\in A_i.
\ea\right.
\eea
Then we may define similarly the raw set valued Hamiltonian and set valued Hamiltonian:
\bea
\label{setH}
\left.\ba{c}
\dis \dbH_0(\th) := \big\{ h(\th, a^*): a^*\in \cE(\th)\big\}\subset \dbR^N;\q \dbH(\th) := \bigcap_{\e>0}\dbH_\e(\th)\subset \dbR^N,\\
\dis \mbox{where}\q \dbH_\e(\th) := \Big\{y\in \dbR^N: |y-h(\th, a^\e)|<\e~\mbox{for some}~a^\e\in\cE_\e(\th)\Big\}. 
\ea\right.
\eea

\ms
\begin{rem}
\label{rem-barh}
Introduce, for any $\th\in\dbX$, $a\in A$, $\a\in \cA$, and $i=1,\cds, N$,
\bea
\label{olhY}
\left.\ba{c}
\dis \ol h_i(\th^i, a) := \sup_{\tilde a_i\in A_i} h_i\big(\th^i, (a^{-i}, \tilde a_i)\big);\\
\dis  \ol Y^{\a,i}_t = g_i(X_T) + \int_t^T \ol h_i(s, X_s, \ol Z^{\a,i}_s, \a_s) ds - \int_t^T \ol Z^{\a,i}_s \cd dB_s.
\ea\right.
\eea
Note that $\ol h_i$ depends only on $a^{-i}$, not on $a_i$.
Under Assumption \ref{assum-standing}, $\ol h$ inherits the regularity \reff{hreg} of $h$ and in particular is  uniformly Lipschitz continuous in $z_i$. By standard BSDE theory (c.f. \cite{Zhangbook}), we have $\ol Y^{\a,i}_0 = \sup_{\tilde \a^i\in \cA_i}J_i(\a^{-i}, \tilde\a^i)$.  Then
\bea
\label{olYa*}
\left.\ba{c}
\a^*\in \cE \q\mbox{if and only if}\q Y^{\a^*,i}_0 = \ol Y^{\a^*, i}_0,~ \forall i=1,\cds, N;\ss\\
 \a^\e\in \cE_\e \q\mbox{if and only if}\q Y^{\a^\e,i}_0 \ge \ol Y^{\a^\e, i}_0 -\e,~ \forall i=1,\cds, N;\ss\\
 a^* \in \cE(\th) \q\mbox{if and only if}\q h_i(\th^i, a^*)  = \ol h_i(\th^i, a^{*}),\q \forall i=1,\cds, N;\ss\\
 a^\e \in \cE_\e(\th) \q\mbox{if and only if}\q h_i(\th^i, a^\e) \ge \ol h_i(\th^i, a^{\e})-\e,\q \forall i=1,\cds, N.
 \ea\right.
 \eea
\end{rem}

Denote by cl the closure of a set. By \cite[Proposition 4]{feinstein2020dynamic} and \reff{YZa}, we have

\begin{prop}
\label{prop-VHcompact} Under Assumption \ref{assum-standing}, for arbitrary $\th\in\dbX$, we have 

\no(i) For any $\e>0$, $\dbV_\e$ and $\dbH_\e(\th)$ are open.

\no(ii) For any $0\le \e_1 < \e_2$, $cl(\dbV_{\e_1}) \subset \dbV_{\e_2}$, $cl(\dbH_{\e_1}(\th)) \subset \dbH_{\e_2}(\th)$.

\no(iii) $\dbV$ and $\dbH(\th)$ are compact, and $|\eta|\le C[1+|z|]$ for all $\eta\in \dbH(\th)$.
\end{prop}
We note that the linear growth of $\dbH(\th)$ is due to the boundedness assumption of $b, f$.  Another typical case in the stochastic control and game literature is that $f_i$ is strictly concave in $a_i$ (and thus unbounded). In this case $\dbH$ may have quadratic growth in $z$ and, as  in stochastic control literature, will become a lot more challenging technically.  

The following two results won't be used in this paper, but are interesting in their own rights.   We postpone their proofs to Appendix.

\begin{prop}\label{prop-HH0}
Under Assumption \ref{assum-standing}, if $A$ is compact, then $\dbH(\th)=\dbH_0(\th)$ for all $\th\in \dbX$.
\end{prop}

For a process $\eta\in \dbL^2(\dbF; \dbR^N)$, let $(Y^\eta, Z^\eta)$ denote the solution to the following BSDE:
\bea
\label{BSDEeta}
Y^{\eta,i}_t = g_i(X_T) + \int_t^T \eta^i_s ds - \int_t^T Z^{\eta,i}_s \cd dB_s, \q i=1,\cds, N.
\eea

\begin{prop}
\label{prop-VH0}
Under Assumption \ref{assum-standing}, $y\in \dbV_0$ if and only if $y = Y^\eta_0$ for some $\eta\in  \dbL^2(\dbF; \dbR^N)$ such that 
\bea
\label{etaH0}
\eta_t \in \dbH_0(t, X_t, Z^\eta_t),~ dt\times d\dbP\mbox{-a.s.}
\eea
\end{prop}
We remark that the if direction has been proved  by \cite{zbMATH01066318} in a different form. 

\begin{rem}
\label{rem-VH0}
We may state Proposition \ref{prop-VH0} alternatively as follows. 
Let $\cH_0$ denote the set of progressively measurable functions $H: [0, T]\times C([0, T];\dbR^d)\times \dbR^{dN} \to \dbR^N$ such that $H$ is   adapted in $\bx\in C([0, T];\dbR^d)$ and $H(t, \bx,z)\in \dbH_0(t, \bx_t, z)$ for all $(t, \bx,z)$.  For each $H\in \cH_0$, let $\dbY^H_0$ denote the set of $Y^H_0$ where $(Y^H, Z^H)$ solves the following path dependent BSDE:
\bea
\label{BSDEH}
Y^{H,i}_t = g_i(X_T) + \int_t^T H_i(s, X_\cd, Z^{H,i}_s) ds - \int_t^T Z^{H,i}_s \cd dB_s,\q i=1,\cds, N.
\eea
We remark that neither the existence nor the uniqueness of the above BSDE is guaranteed, and thus $\dbY^H_0$ could be empty or have multiple values. Then Proposition \ref{prop-VH0} is equivalent to:
\bea
\label{VH0}
V_0 = \bigcup_{H\in \cH_0} \dbY^H_0.
\eea

We emphasize that in general $H$ has to be path dependent in $X$, as we will explain in Remark \ref{rem-path} below. 
We prefer the form \reff{BSDEeta}-\reff{etaH0} than \reff{BSDEH} because \reff{BSDEeta} is always wellposed. However, in \S\ref{sect-Lipschitz} below, we will restrict to Lipschitz continuous $H$ and thus \reff{BSDEH} is also wellposed, and in that case it is more convenient to use \reff{BSDEH} as we will indeed do.
\end{rem} 

Our goal of this paper is to extend the above relationship to $\dbV$ and $\dbH$. For this purpose, we introduce the Hausdorff distance $d(\dbH_1, \dbH_2)$  for two sets $\dbH_1, \dbH_2$:
\bea
\label{dbHdistance}
d(\dbH_1,\dbH_2)=\big[\sup_{y_1\in \dbH_1} d(y_1, \dbH_2)\big] \vee \big[\sup_{y_2\in \dbH_2} d(y_2, \dbH_1)\big],\q d(y_1, \dbH_2) := \inf_{y_2\in \dbH_2}|y_1-y_2|.
\eea
Throughout the paper, we shall also impose the following assumptions.

\begin{assum}
\label{assum-dbH}
 $\dbH(\th) \neq \emptyset$ for all $\th\in\dbX$, and $\dbH$ is continuous in $\th$ under $d$.
\end{assum}

\begin{rem}
\label{rem-dbH}
(i) Since we allow $\dbV$ to be empty, the assumption $\dbH(\th) \neq \emptyset$ is redundant. However, in the case $\dbH(\th) =\emptyset$ some results in the paper need to be stated more carefully. Since the nonempty case is more interesting in practice, we impose this assumption for the convenience of presentation. We shall emphasize that it is a lot easier to verify $\dbH(\th) \neq \emptyset$ then to verify $\dbV \neq \emptyset$.

\ss
\no(ii) The continuity of $\dbH$ is not a trivial requirement. It will be very interesting to explore the possible relaxation of this assumption. 
\end{rem}

 In the next section we assume $\dbH$ is continuous in all variables, while in \S\ref{sect-Lipschitz} we assume further that $\dbH$ is uniformly Lipschitz continuous in $z$.

\section{The case with continuous   $\dbH$}
\label{sect-continuous}

We first establish the following proposition which will be crucial for our main results. For any compact subset $K\subset\subset \dbX$ and $\e>0$, denote
      \bea
      \label{rhoKe}
        \rho_K(\e):=\sup_{\th\in K}d\big(\dbH_\e(\th),\dbH(\th)\big).
    \eea

\begin{prop}
\label{prop-Hcont}
Under Assumption \ref{assum-standing}, for any compact subset $K\subset\subset \dbX$, $\dis\lim_{\e\to 0}\rho_K(\e)  =0$   if and only if $\dbH$ is (uniformly) continuous  on $K$.
\end{prop}
Note that $\dbH_\e(\th)$ is monotone in $\e$, so the if direction can be viewed as a set valued version of  Dini's lemma. 
We shall postpone its proof to \S\ref{sect-Hcont} below. The next result follows from somewhat standard arguments, and we shall postpone its proof to Appendix.

\begin{prop}
\label{prop-construction}
Let Assumptions \ref{assum-standing} and \ref{assum-dbH} hold.  Then, for any $\eta\in \dbL^2(\dbF; \dbR^d)$, $Z\in \dbL^2(\dbF; \dbR^{dN})$, and any $\e>0$, there exist $\a^\e\in \cA$ and $Z^\e\in \dbL^2(\dbF; \dbR^{dN})$ such that
\bea
\label{construction}
\left.\ba{c}
\dis \a^\e_t \in \cE_\e(t, X_t, Z^\e_t),\q  \dbE\Big[\int_0^T |Z^\e_t-Z_t|^2dt\Big] \le \e,\ss\\
\dis \dbE\Big[\int_0^T \Big|\big|\eta_t - h(\Th_t, \a^\e_t)\big|-d\big(\eta_t, \dbH(\Th_t)\big)\Big|^2dt\Big] \le \e,\q\mbox{where}\q \Th_t := (t, X_t, Z_t).
\ea\right.
\eea
\end{prop}

The main result of this section is as follows.  
\begin{thm}
\label{thm-general}
Let Assumptions \ref{assum-standing} and \ref{assum-dbH} hold.   Then $\dis \dbV = \bigcap_{\e>0} \big\{Y^\eta_0: \eta\in \Xi_\e\big\}$, where
\bea
\label{Xie}
\Xi_\e := \Big\{\eta\in \dbL^2(\dbF; \dbR^N): \dbE\Big[\int_0^T d^{3\over 2}\big(\eta_t, \dbH(t, X_t, Z^\eta_t)\big) dt\Big] \le \e\Big\}.
\eea
\end{thm}

\begin{rem}
\label{rem-general}
(i) The order ${3\over 2}$ here is for technical reasons, which can be replaced with any order between $1$ and $2$ (exclusively). 

\ss
\no(ii) When $\e=0$, $\eta\in \Xi_0$ means $\eta_t\in \dbH(t, X_t, Z^\eta_t)$. It will be interesting to explore under what conditions $\dbV = cl\big\{Y^\eta_0: \eta\in \Xi_0\big\}$. In particular, we shall provide a positive answer in the next section under additional conditions on  $\dbH$. We remark that, in light of Propositions \ref{prop-HH0} and \ref{prop-VH0}, when $A$ is compact this claim is equivalent to $\dbV = cl(\dbV_0)$. 

\ss
\no (iii) It should be noted that we are using the value $Y^\eta_0$ rather than its $\e$-neighborhood. In other words, for $y\in \dbV$, for any $\e>0$, there exists $\eta^\e\in \Xi_\e$ such that $Y^{\eta^\e}_0=y$. 
\end{rem}

\no{\bf Proof of Theorem \ref{thm-general}.}  Denote $\dis\tilde V := \bigcap_{\e>0} \big\{Y^\eta_0: \eta\in \Xi_\e\big\}$.

(i) We first prove $\tilde V \subset \dbV$. Fix $y\in \tilde \dbV$ and $\e>0$. Then $y = Y^\eta_0$ for some $\eta\in \Xi_\e$. By Proposition \ref{prop-construction} there exist $\a^\e$ and $Z^\e$ such that \reff{construction} holds with $Z=Z^\eta$. 
Denote $\Th_t := (t, X_t, Z^\eta_t)$, $\Th^\e_t := (t, X_t, Z^\e_t)$, $\D Y_t := Y^{\eta}_t-Y^{\a^\e}_t$, $\D Z_t := Z^{\eta}_t- Z^{\a^\e}_t$. Then we have
\beaa
\D Y^i_t = \int_t^T \Big[\eta^i_t - h_i(\Th^i_t, \a^{\e}_t)+ \b^i_s \D Z^i_s\Big] ds- \int_t^T \D Z^i_s \cd dB_s,
\eeaa
where $\b^i$ is bounded due to the Lipschitz continuity of $h$ in $z$ as in \reff{hreg}. Denote 
\bea
\label{Mbeta}
M^{\b^i}_T := \exp\Big(\int_0^T \b^i_s \cd dB_s - {1\over 2}\int_0^T |\b^i_s|^2 ds\Big). 
\eea
Then $\dbE\big[|M^{\b^i}_T|^p + |M^{\b^i}_T|^{-p}\big]\le C_p <\infty$, for any $p\ge 1$. Thus, by \reff{construction} and \reff{Xie},
\bea
\label{DYi0}
&&\dis\!\!\!\!\!\!\!\!\!\!\!\! |\D Y^i_0| = \Big|\dbE\Big[M^{\b^i}_T \int_0^T \big[\eta^i_t - h_i(\Th^i_t, \a^{\e}_t)\big]dt\Big]\Big|\le \dbE\Big[M^{\b^i}_T \int_0^T \big|\eta_t - h(\Th_t, \a^{\e}_t)\big|dt\Big]\nonumber\\
&&\dis \!\!\!\!\!\!\!\!\!\!\!\!\le \dbE\Big[M^{\b^i}_T \int_0^T \Big[\big||\eta_t - h(\Th_t, \a^{\e}_t)| -  d(\eta_t, \dbH(\Th_t))\big| + d(\eta_t, \dbH(\Th_t))\big| \Big] dt\Big] \nonumber\\
&&\dis\!\!\!\!\!\!\!\!\!\!\!\! \le C \Big(\dbE\Big[\int_0^T \big||\eta_t - h(\Th_t, \a^{\e}_t)| -  d(\eta_t, \dbH(\Th_t))\big|^2dt\Big]\Big)^{1\over 2} + C\Big(\dbE\Big[\int_0^Td^{3\over 2}(\eta_t, \dbH(\Th_t))\big| dt\Big]\Big)^{2\over 3} \nonumber\\
&&\dis\!\!\!\!\!\!\!\!\!\!\!\! \le C[\e^{1\over 2} + \e^{2\over 3}] \le C\sqrt{\e}. 
\eea
Here we assume without loss of generality that $\e\le 1$.
Moreover, since $\a^\e_t \in \cE_\e(\Th^\e_t)$, by \reff{olYa*} we have $0\le \ol h^i(\Th^{\e, i}_t, \a^{\e}_t)- h^i(\Th^{\e,i}_t, \a^{\e}_t)\le \e$. Thus, \reff{hreg} leads to
\beaa
 \big|\eta^i_t - \ol h^i(\Th^i_t, \a^{\e}_t)\big| \le \big|\eta^i_t - h^i(\Th^{i}_t, \a^{\e}_t)\big| + C|Z^\e_t-Z_t|+ \e.
\eeaa
Then by \reff{construction}  similarly we have $\big|Y^{\eta, i}_0-\ol Y^{\a^\e,i}_0\big| \le  C \sqrt{\e}$. This, together with \reff{DYi0}, implies $ \big|Y^{\a^\e, i}_0-\ol Y^{\a^\e,i}_0\big| \le  C \sqrt{\e}$. Then,
by \reff{olYa*} we see that $\a^\e\in \cE_{C \sqrt{\e}}$ and $|y-Y^{\a^\e}_0| = |Y^\eta_0-Y^{\a^\e}_0|\le C \sqrt{\e}$. That is, $y\in \dbV_{C \sqrt{\e}}$. Since $\e>0$ is arbitrary, we obtain $y\in \dbV$.

(ii) We next show that $\dbV \subset \tilde \dbV$. Fix $y\in \dbV$ and $\e>0$. Then $|y-Y^{\a^\e}_0| \le \e$ for some $\a^\e\in \cE_\e$.  Set $\Th_t := (t, X_t, Z^{\a^\e}_t)$, $\eta_t := h(\Th_t, \a^\e_t)$. Then it is clear that $Y^{\a^\e} = Y^\eta$, $Z^{\a^\e} = Z^\eta$.

Denote $\D \ol Y_t := \ol Y^{\a^\e}_t - Y^{\a^\e}_t$, $\D \ol Z_t := \ol Z^{\a^\e}_t - Z^{\a^\e}_t$. Again by \reff{hreg} we have
\beaa
\D \ol Y^i_t = \int_t^T \Big[\ol h_i(\Th^i_s, \a^{\e}_s) - h_i(\Th^i_s, \a^\e_s) + \b^i_s \D \ol Z^i_s\Big] ds- \int_t^T \D\ol Z^i_s \cd dB_s,
\eeaa
for some appropriate bounded process $\b^i$. Recall \reff{Mbeta}.  By \reff{olYa*} we have, 
\bea
\label{DY0<e}
\e \ge \D \ol Y^i_0 = \dbE\Big[M^{\b^i}_T \int_0^T \D \eta^i_s ds\Big],~\mbox{where}~  \ol\eta^i_s:=\ol h_i(\Th^i_s, \a^{\e}_s),~ \D \eta^i_s:=\ol\eta^i_s  - \eta^i_s\ge 0.
\eea
 For each $R>0$, recall \reff{rhoKe} and denote
\beaa
&\dis  E^\e_t := \bigcap_{i=1}^N \big\{\D \eta^i_t\le \sqrt{\e}\big\};\\
&\dis K_R:= \big\{(x, z): |x|+ |z|\le R\big\},\q  \rho_R(\e) := \rho_{[0,T]\times K_R}(\sqrt{\e}),\q E^R_t := \big\{(X_t, Z^\eta_t) \in K_R\big\}.
\eeaa
Then, by \reff{olYa*}, $\a^\e_t\in \cE_{\sqrt{\e}}(\Th_t)$ and  thus $\eta_t \in \dbH_{\sqrt{\e}}(\Th_t)$ on $E^\e_t$. Therefore, 
\beaa
d(\eta_t, \dbH(\Th_t)) \le d\Big(\dbH_{\sqrt{\e}}(\Th_t), \dbH(\Th_t)\Big)  \le \rho_R(\e),\q\mbox{on}~ E^\e_t \cap E^R_t.
\eeaa

Note that $|h|\le C(1+|z|)$, $|g|\le C$. Applying standard BSDE estimates on \reff{YZa} we get
\beaa
\dbE\Big[\int_0^T [|Z^\eta_t|^2 + |\eta_t|^2] dt\Big] \le C,\q d(\eta_t, \dbH(\Th_t)) \le C[|\eta_t|+|Z^\eta_t|].
\eeaa
Then
\beaa
&&\dis \dbE\Big[\int_0^T d^{3\over 2}(\eta_t, \dbH(\Th_t)) dt\Big] \le T \rho_R^{3\over 2}(\e) + C\dbE\Big[\int_0^T \big[|\eta_t|^{3\over 2}+|Z^\eta_t|^{3\over 2}\big]\1_{(E^\e_t)^c \cup (E^R_t)^c} dt\Big]\\
&&\dis \le C \rho_R^{3\over 2}(\e) + C\dbE\Big[\int_0^T \big[|\eta_t|^{3\over 2}+|Z^\eta_t|^{3\over 2}\big] \big[\sum_{i=1}^N \1_{\{\D \eta^i_t>\sqrt{\e}\}} + \1_{\{|X_t|+|Z^\eta_t|>R\}}\big] dt\Big]
\eeaa
Note that
\beaa
&&\dis \dbE\Big[\int_0^T \big[|\eta_t|^{3\over 2}+|Z^\eta_t|^{3\over 2}\big]  \1_{\{|X_t|+|Z^\eta_t|>R\}} dt\Big] \le {1\over \sqrt{R}}\dbE\Big[\int_0^T \big[|\eta_t|^{3\over 2}+|Z^\eta_t|^{3\over 2}\big]  \big[|X_t|^{1\over 2}+|Z^\eta_t|^{1\over 2}\big] dt\Big] \\
&&\dis\qq \le {C\over \sqrt{R}}\dbE\Big[\int_0^T \big[|\eta_t|^2+|Z^\eta_t|^2+|X_t|^2\big] dt\Big] \le {C\over \sqrt{R}};\\
&&\dis \dbE\Big[\int_0^T \big[|\eta_t|^{3\over 2}+|Z^\eta_t|^{3\over 2}\big] \1_{\{\D \eta^i_s>\sqrt{\e}\}} dt\Big]\le  \e^{-{1\over 16}}\dbE\Big[\int_0^T \big[|\eta_t|^{3\over 2}+|Z^\eta_t|^{3\over 2}\big] (\D \eta^i_s)^{1\over 8} dt\Big]\\
&&\dis\qq = \e^{-{1\over 16}}\dbE\Big[(M^{\b^i}_T)^{-{1\over 8}}\int_0^T \big[|\eta_t|^{3\over 2}+|Z^\eta_t|^{3\over 2}\big] (M^{\beta^i}_T\D \eta^i_t)^{1\over 8} dt\Big]\\
&&\dis\qq \le \e^{-{1\over 16}}\Big(\dbE\big[(M^{\b^i}_T)^{-1}\big]\Big)^{1\over 8}\Big(\dbE\Big[\int_0^T\big[|\eta_t|^2+|Z^\eta_t|^2\big]dt\Big]\Big)^{3\over 4} \Big(\dbE\Big[M^{\b^i}_T\int_0^T\D \eta^i_t dt\Big]\Big)^{1\over 8} \\
&&\dis\qq \le C\e^{-{1\over 16}}\Big(\dbE\Big[M^{\b^i}_T\int_0^T\D \eta^i_t dt\Big]\Big)^{1\over 8}\le C\e^{1\over 16},
\eeaa
where the last inequality is due to \reff{DY0<e}. Put together, we have 
\bea
\label{detaHTh}
 \dbE\Big[\int_0^T d^{3\over 2}(\eta_t, \dbH(\Th_t)) dt\Big] \le C \big[\rho_R^{3\over 2}(\e) + \e^{1\over 16} + R^{-{1\over 2}}\big].
\eea

Introduce further $\tilde \eta_t := \eta_t + {1\over T}(y- Y^{\a^\e}_0)$. Then one can easily see that 
\beaa
Z^{\tilde \eta} = Z^\eta, \q Y^{\tilde \eta}_t = Y^{\eta}_t + {T-t\over T}(y- Y^{\a^\e}_0), 
\eeaa
and in particular, $Y^{\tilde \eta}_0 = Y^{\eta}_0 + (y- Y^{\a^\e}_0) = y$. Moreover, by \reff{detaHTh} and recalling $|y-Y^{\a^\e}_0| \le \e$,
\beaa
 \dbE\Big[\int_0^T d^{3\over 2}(\tilde \eta_t, \dbH(t, X_t, Z^{\tilde \eta}_t)) dt\Big] &\le&  C \dbE\Big[\int_0^T \big[d^{3\over 2}( \eta_t, \dbH(\Th_t)) + |\tilde \eta_t-\eta_t|^{3\over 2}\big]dt\Big] \\
 &\le& C_0 \big[\rho_R^{3\over 2}(\e) + \e^{1\over 16} + R^{-{1\over 2}} + \e^{3\over 2}\big],
\eeaa
for some  constant $C_0$.
Now for any $\e_0>0$, by first setting $R := ({2C_0\over \e_0})^2$ so that $C_0 R^{-{1\over 2}} \le {1\over 2}\e_0$, and then choosing $\e>0$ small enough such that $C_0 \big[\rho_R^{3\over 2}(\e) + \e^{1\over 16} + \e^{3\over 2}\big] \le {1\over 2}\e_0$, we have 
\beaa
 \dbE\Big[\int_0^T d^{3\over 2}(\tilde \eta_t, \dbH(t, X_t, Z^{\tilde \eta}_t)) dt\Big] \le \e_0.
\eeaa
That is, $\tilde \eta\in \Xi_{\e_0}$ and $y= Y^{\tilde \eta}_0$. By the arbitrariness of $\e_0>0$, we see that $y\in \tilde \dbV$.
\qed

\subsection{Proof of Proposition \ref{prop-Hcont}}
\label{sect-Hcont}
First, for any fixed $\th\in\dbX$, by Proposition \ref{prop-VHcompact} we have $\dis\dbH(\th) = \bigcap_{\e>0} cl(\dbH_\e(\th))$ and each $cl(\dbH_\e(\th))$ is compact. Then it's clear that 
\bea
\label{limitHe}
\lim_{\e\to 0} d(\dbH_\e(\th),\dbH(\th))=\lim_{\e\to 0} d\big(cl(\ol\dbH_\e(\th)), \dbH(\th)\big) =0.
\eea

{\bf Step 1.} We first prove the ``if'' direction. Assume by contradiction that $\dbH$ is continuous on $K$ but $\dis c_0:= \limsup_{\e\to 0}\rho_K(\e)  >0$. Note that $\dbH(\th) \subset \dbH_\e(\th)$ and thus $d\big(y, \dbH_\e(\th)\big)=0$ for all $y\in \dbH(\th)$ and $\e<0$. Then there exist $\e_n \downarrow 0$, $\th_n \in K$, and $y_n \in \dbH_{\e_n}(\th_n)$ such that $d(y_n, \dbH(\th_n)) \ge {c_0\over 2}$. Since $K$ is compact, by otherwise considering a subsequence, we may assume $\th_n \to \th_*\in K$. By the continuity of $\dbH$, for $n$ large enough, we have 
\bea
\label{dHn}
d(\dbH(\th_n),\dbH(\th_*)) \le {c_0\over 4},\q\mbox{and then}\q  d(y_n, \dbH(\th_*)) \ge {c_0\over 4}.
\eea 
Since $y_n \in \dbH_{\e_n}(\th_n)$, there exists $a^n \in \cE_{\e_n}(\th_n)$ with $|y_n - h(\th_n, a^n)|< \e_n$. By \reff{olYa*} we have
\beaa
 h_i(\th^i_n, a^n) \ge \ol h_i(\th^i_n, a^{n})-\e_n,\q  i=1,\cds, N.
\eeaa
By \reff{YZa} and \reff{olhY}, since $K$ is compact, we see that $h, \ol h$ are uniformly continuous in $\th$ for $\th\in K$, with a modulus of continuity function $\rho$ which may possibly depend on $K$.  Then
\beaa
h_i(\th^i_*, a^n) \ge \ol h_i(\th^i_*, a^{n})-\e_n- 2\rho(|\th_n-\th_*|),\q i=1,\cds, N.
\eeaa
That is, $a^n \in \cE_{\e_n+2\rho(|\th_n-\th_*|)}(\th_*)$. Moreover, note that 
\beaa
|y_n-h(\th_*, \a^n)| < \e_n + |h(\th_n, \a^n)-h(\th_*, \a^n)| \le \e_n + \rho(|\th_n-\th_*|).
\eeaa
Then $y_n \in \dbH_{\e_n+2\rho(|\th_n-\th_*|)}(\th_*)$, and thus, by \reff{limitHe},
\beaa
d(y_n, \dbH(\th_*)) \le d\big( \dbH_{\e_n+2\rho(|\th_n-\th_*|)}(\th_*), \dbH(\th_*)) \to 0,\q\mbox{as}\q n\to\infty. 
\eeaa
This contradicts with \reff{dHn}, and thus $\dis\lim_{\e\to 0}\rho_K(\e)  =0$.

\ms
{\bf Step 2.} We next prove the ``only if'' direction. Assume $\dis\lim_{\e\to 0}\rho_K(\e)  =0$ and let $\th_n, \th_*\in K$ be such that $\dis \lim_{n\to\infty}|\th_n - \th_*|\to 0$. Fix $\e>0$. First, for any $y\in \dbH(\th_*)\subset  \dbH_{\e\over 2}(\th_*)$, there exists $a^\e\in \cE_{\e\over 2}(\th_*)$ such that $|y- h(\th_*, a^\e)|< {\e\over 2}$.  For each $i=1,\cds, N$, by  \reff{olYa*} we have  $h_i(\th^i_*, a^\e) \ge \ol h_i(\th^i_*, a^{\e})-{\e\over 2}$. By the uniform continuity \reff{hreg} of $h, \ol h$ in $\th\in K$, there exists $n_\e$ independent of $a^\e$ such that, for all $n\ge n_\e$,
\beaa
h_i(\th^i_n, a^\e) \ge \ol h_i(\th^i_n, a^{\e})-\e,\q\mbox{and}\q |h(\th_*, \a^\e)-h(\th_n, \a^\e)| \le {\e\over 2}.
\eeaa
Then $a^\e\in \cE_\e(\th_n)$, and
\beaa
|y-h(\th_n, \a^\e)| < {\e\over 2} + |h(\th_*, \a^\e)-h(\th_n, \a^\e)| \le \e.
\eeaa
This implies $y\in \dbH_\e(\th_n)$, and thus $d(y, \dbH(\th_n)) \le d(\dbH_\e(\th_n), \dbH(\th_n)) \le \rho_K(\e)$, for all $n\ge n_\e$.  Therefore, $\sup_{y\in \dbH(\th_*)} d(y, \dbH(\th_n)) \le \rho_K(\e)$, for all $n\ge n_\e$.

On the other hand, by the same arguments we have  $\dis \sup_{y_n\in \dbH(\th_n)} d(y, \dbH(\th_*)) \le \rho_K(\e)$ for the same $n_\e$ and for all $n\ge n_\e$.
Then $d(\dbH(\th_n), \dbH(\th_*)) \le \rho_K(\e)$, for all $n\ge n_\e$. This implies that $\dis\limsup_{n\to\infty} d(\dbH(\th_n), \dbH(\th_*)) \le \rho_K(\e)$. Sending $\e\to 0$, we obtain $\dis\lim_{n\to\infty} d(\dbH(\th_n), \dbH(\th_*))=0$.
\qed

\section{The case with Lipschitz continuous  $\dbH$}
\label{sect-Lipschitz}

For $L\ge 0$, let $\cH_L$ denote the set of progressively measurable functions $H: [0, T]\times C([0, T]; \dbR^d) \times \dbR^{dN}\to \dbR^N$ such that $H$ is adapted in $\bx\in C([0, T]; \dbR^d)$, and

\ss
$\bullet$ $H(t, \bx_{\cd\wedge t}, z) \in \dbH(t, \bx_t, z)$, for all $(t, \bx, z)$;

\ss
$\bullet$ $H$ is uniformly Lipschitz continuous in $z$ with Lipschitz constant $L$. 

\ss
\no Note that we do not require $H$ to be continuous in $(t, \bx)$. By Proposition \ref{prop-VHcompact} (iii),  $\dbH(t, x, 0)$ is uniformly bounded, then \reff{BSDEH} is wellposed for any $H\in \cH_L$.  In this section we first establish our main result under the following assumption. The essence of this assumption is that $\dbH$ is Lipschitz continuous in $z$, as we will discuss more  in \S\ref{sect-LipH} below.  

\begin{assum}
\label{assum-Lipschitz}
 There exists $L\ge 0$ such that, for any $\eta\in \dbL^2(\dbF; \dbR^N)$, $Z\in \dbL^2(\dbF; \dbR^{dN})$, and any $\e>0$, there exists $H^\e \in \cH_L$ such that 
 \bea
 \label{construction2}
 \dbE\Big[\int_0^T\Big| \big|\eta_t - H^\e(t, X_\cd, Z_t)\big| - d\big(\eta_t, \dbH(t, X_t, Z_t)\big)\Big|^2dt\Big] \le \e.
 \eea
\end{assum}
\begin{thm}
\label{thm-Lipschitz}
Let Assumptions \ref{assum-standing}, \ref{assum-dbH}, and \ref{assum-Lipschitz} hold true  with constant $L$. Then
\bea
\label{VYH}
\dbV = {\rm cl}\big\{Y^H_0: H \in \cH_L\big\}.
\eea
\end{thm}
\proof We first prove $\supset$, and we emphasize that this inclusion actually does not require Assumption \ref{assum-Lipschitz}.
 Fix $H \in \cH_L$, and set  $\eta_t:= H(t, X_\cd, Z^H_t)\in \dbH(t, X_t, Z^H_t)$. By  \reff{BSDEeta} and \reff{BSDEH}, it is clear that $Y^\eta = Y^H$, $Z^\eta = Z^H$. In particular, this implies that $d\big(\eta_t, \dbH(t, X_t, Z^\eta_t)\big)=0$, and thus $\eta\in \Xi_\e$ for all $\e>0$. Then it follows from Theorem \ref{thm-general} that $Y^H_0=Y^\eta_0 \in \dbV$.

 We next prove $\subset$. Fix $y\in \dbV$ and $\e >0$. By Theorem \ref{thm-general}, $y= Y^\eta_0$ for some $\eta \in \Xi_\e$. Let $H^\e\in \cH_L$ be as in Assumption \ref{assum-Lipschitz} with $Z=Z^\eta$. Combining \reff{Xie} and \reff{construction2} we have
 \beaa
 \dbE\Big[\int_0^T |\D \eta^\e_t|^{3\over 2} dt\Big] \le C\e^{3\over 4},\q\mbox{where}\q \D \eta^\e_t:= \eta_t - H^\e(t, X_\cd, Z^\eta_t).
 \eeaa
Now denote $\D Y^\e :=  Y^\eta - Y^{H^\e}$, $\D Z^\e := Z^\eta - Z^{H^\e}$. Since $H^\e\in \cH_L$, we have
\beaa
\D Y^i_t = \int_t^T \big[\D \eta^i_s + \beta^i_s \D Z^i_s\big]ds - \int_t^T \D Z^i_s dB_s,
\eeaa
where $|\beta^i|\le L$. Thus, recalling \reff{Mbeta},
\bea
\label{yYHe}
|y- Y^{H^\e}_0| = |\D Y^i_0| = \Big|\dbE\Big[M^{\b^i_T} \int_0^T \D \eta^i_s ds\Big] \Big|\le C \Big(\dbE\Big[ \int_0^T |\D \eta^i_s|^{3\over 2} ds\Big] \Big)^{2\over 3} \le C\sqrt{\e}.
\eea
 By the arbitrariness of $\e>0$, we see that $y \in cl\{Y^H_0: H\in \cH_L\}$.
\qed

\begin{rem}
\label{rem-path}
In the proof of $\subset$ at above, given $\eta \in \Xi_\e$, in general, $\eta_t$ may not be a deterministic function of $(X_t, Z^{\a^\e}_t)$. That is, if is possible that $(X_t, Z^{\a^\e}_t)(\o^1) = (X_t, Z^{\a^\e}_t)(\o^2)$ but $\eta_t(\o^1) \neq \eta_t(\o^2)$ for different $\o^1, \o^2\in \O$.  Thus we have to allow $H^\e$ to depend on the paths of $X$. This is also consistent with \cite[Corollary 1 and Proposition 3]{feinstein2020dynamic}  that in general the DPP for the dynamic set value function $\dbV(t,x)$ has to allow path dependent approximate equilibria.  We also refer to Remark \ref{rem-path2} for further remarks on this issue.
\end{rem}

In the rest of this section we present a few important cases where Assumption \ref{assum-Lipschitz} can be verified and hence $\dbV$ can be characterized by Theorem \ref{thm-Lipschitz}. 

\subsection{The singleton case}
First, when $\dbH(t,x,z) = \{H(t,x,z)\}$ is a singleton, it is straightforward to verify Assumption \ref{assum-Lipschitz} with $H^\e\equiv H$.  Then the following result  is a direct consequence of Theorem \ref{thm-Lipschitz} and the well known connection between BSDEs and PDEs.
\begin{thm}
\label{thm-singleton}
Let Assumptions \ref{assum-standing} and \ref{assum-dbH} hold. Assume $\dbH(t,x,z) = \{H(t,x,z)\}$ and $H$ is  uniformly Lipschitz continuous in $z$. Then $\dbV = \{Y^H_0\}$ is also a singleton.

Moreover, $Y^H_0 = u(0, 0)$, where $u$ solves the following system of parabolic PDEs:
\bea
\label{PDE}
\pa_t u^i + {1\over 2}\tr(\pa_{xx} u^i) + H_i(t, x, \pa_x u)=0, \q u^i(T,x) = g_i(x),\q i=1,\cds, N.
\eea
\end{thm} 
In this case it is natural to call $H$ the Hamiltonian  and $Y^H_0$ the game value of the game.
\begin{rem}
\label{rem-singleton}
 Theorem \ref{thm-singleton} is a strict generalization of the standard stochastic control problems and two person zero sum game problems under Isaacs conditions. In particular, as in those two cases, we may study the game value $Y^H_0$ without the existence of Nash equilibrium.

(i) When $N=1$, the problem reduces to a standard stochastic control problem. It is obvious that $\dbH(t,x,z) = \{H(t,x,z)\}$, where $H(t,x,z) := \sup_{a\in A} h(t,x,z,a)$, and $\dbV = \{u(0,0)\}$, where $u$ is the solution to the standard HJB equation \reff{PDE} with $N=1$. 

(ii) When $N=2$ and $(f_2, g_2) = - (f_1, g_1)$, the problem reduces to a two person zero sum game. There are two cases.

{\it Case 1.} The Isaacs condition holds: 
\bea
\label{Isaacs}
\inf_{a_2\in A_2} \sup_{a_1\in A_1} h_1(t,x,z_1, a_1, a_2) = \sup_{a_1\in A_1}\inf_{a_2\in A_2}  h_1(t,x,z_1, a_1, a_2).
\eea
Denote the above as $H_1(t,x,z_1)$.  Then $\dbH(t,x,z) = \big\{(H_1(t,x,z_1), H_2(t,x,z_2))\big\}$ and $\dbV = \big\{\big(u_1(0,0), u_2(0,0)\big)\big\}$ are singletons with $H_2(t,x, z_2) = -H_1(t,x,z_2)$, $u_2 = - u_1$,  where
\beaa
 \pa_t u_1 + {1\over 2} \pa_{xx} u_1 + H_1(t,x,\pa_x u_1) =0,\q u_1(T,x) = g_1(x).
\eeaa
In particular, the (unique) game value exists in this case.

{\it Case 2.} The Isaacs condition \reff{Isaacs} fails. Then $\dbH(t,x,z) = \emptyset$ for some $(t,x,z)$ and $\dbV = \emptyset$. In particular, the  game value does not exist in this case.
\end{rem} 

\subsection{The separable case}
The next result considers the situation with multiple $H$. The typical case is that $\dbH(\th) = \big\{H^1(\th), \cds, H^n(\th)\big\}$. But we shall go beyond and introduce the separable class.

\begin{defn}
\label{defn-separable}
Given a constant $L\ge 0$, we say $\dbH$ is $L$-separable if  there exist a sequence of functions $H^n: \dbX \to \dbR^N$, $n\ge 1$, such that
\begin{itemize}
\item For each $\th\in \dbX$, $\dbH(\th) = cl\{H^n(\th), n\ge 1\}$;

\item{} For each $n\ge 1$,   $H^n$ is measurable in $(t,x)$ and uniformly Lipschitz continuous in $z$ with a common Lipschitz constant $L$.
 \end{itemize}
\end{defn}
The $L$-separability is obviously closed under union. 
\begin{lem}
\label{lem-union}
Assume $\dbH = cl \big(\bigcup_{n\ge 1} \dbH_n\big)$ and each $\dbH_n$ is $L$-separable with the same $L$. Then $\dbH$ is  $L$-separable.
\end{lem}

  \begin{thm}
\label{thm-multiple}
Let Assumptions \ref{assum-standing} and \ref{assum-dbH} hold, and $\dbH$ be $L$-separable. Then Assumption \ref{assum-Lipschitz} holds true. Moreover, we may restrict to a subclass of $\cH_L$ in \reff{VYH} as follows:
\bea
\label{VYHI}
\dbV = cl\big\{Y^{H^I}_0: I\in \cI\big\},
\eea
where $\cI$ denotes the set of adapted mappings $I: [0, T]\times C([0, T]; \dbR^d) \to \{1, 2,\cds\}$, and,  
\bea
\label{HI}
H^I(t,\bx,z):= H^{I(t,\bx)}(t,\bx_t,z) = \sum_{n=1}^\infty H^n(t, \bx_t, z) \1_{\{I(t,\bx) = n\}}.
\eea
\end{thm} 
\proof We first note that it is obvious that $H^I\in \cH_L$ for all $I\in \cI$. Then by Theorem \ref{thm-Lipschitz}  we have $Y^{H^I}_0 \in \dbV$. We now prove the opposite inclusion. 

Fix $y\in \dbV$, $\e>0$, and $y= Y^\eta_0$ for some $\eta \in \Xi_\e$. Denote $\Th_t := (t,X_t, Z^\eta_t)$. Since $\dbH(t,x,z) = cl\{H^n(t,x,z), n\ge 1\}$, clearly $d\big(\eta_t, \dbH(\Th_t)\big) = \inf_{n\ge 1} |\eta_t - H^n(\Th_t)|$. Introduce
\beaa
\tilde I^\e_t := \inf\Big\{n\ge 1:  |\eta_t - H^n(\Th_t)| \le d\big(\eta_t, \dbH(\Th_t)\big) + {\sqrt{\e\over T}}\Big\}.
\eeaa
Clearly $\tilde I^\e$ is an $\dbF$-progressively measurable process. Since $\dbF=\dbF^B=\dbF^X$, then there exists $I^\e\in \cI$ such that $\tilde I^\e_t = I^\e(t, X_\cd)$. Denote $H^\e:= H^{I^\e}\in \cH_L$. One can easily check that 
\beaa
0\le |\eta_t - H^\e(\Th_t)| - d\big(\eta_t, \dbH(\Th_t)\big) \le {\sqrt{\e\over T}}.
\eeaa
 Then $H^\e$ satisfies all the requirements in Assumption \ref{assum-Lipschitz} and thus Assumption \ref{assum-Lipschitz} holds true. Moreover, it follows from the arguments in the proof  of Theorem \ref{thm-Lipschitz}, especially by \reff{yYHe}, that $|y-Y^{H^{I^\e}}_0|\le C\sqrt{\e}$. Since $\e>0$ is arbitrary,  we obtain $y\in  cl\big\{Y^{H^I}_0: I\in \cI\big\}$. 
\qed

\begin{rem}
\label{rem-path2}
(i) Although all the $H^n$ in Theorem \ref{thm-multiple} are state dependent in $x$, the $\eta$ and  $I^\e$ in the above proof are in general path dependent, then $ H^{I^\e}$ is also path dependent. In Theorem \ref{thm-singleton}, however, we have only one choice: $I\equiv 1$, and thus $H$ becomes state dependent.   

\ms
\no (ii) In light of \reff{PDE}, in the case of path dependent  $H$, we can also characterize $Y^H_0$ through path dependent PDE system. We refer to \cite[Chapter 11]{Zhangbook} for an introduction of the theory of path dependent PDEs.

\ms

\no(iii) Let $\cH_L^{state}$ denote the subset of state dependent $H\in \cH_L$. In general, for a path dependent $H\in \cH_L$, there is no $\tilde H\in \cH^{state}_L$ such that  the processes coincide: $Y^H=Y^{\tilde H}$. However, note that the mapping $H\in \cH_L\mapsto Y^H_0$ is not one to one. It remains a very interesting open problem that whether or not the following equality could hold true:
\bea
\label{state}
cl\big\{Y^H_0: H\in \cH_L\big\} =  cl\big\{Y^H_0: H\in \cH^{state}_L\big\}.
\eea
That is, given $H\in \cH_L$, can we find $H^n\in \cH^{state}_L$ such that $\lim_{n\to \infty} Y^{H^n}_0 = Y^H_0$?
\end{rem}

\begin{rem}
\label{rem-multiple}
Note that it is a lot easier to analyze $\dbH$ than $\dbV$. Under the conditions of Theorem \ref{thm-multiple}, we effectively reduce the  game problem to a control problem with  control $I\in \cI$. This feature is crucial for the work \cite{IseriZhangGamePDE}.
\end{rem}

\subsection{Further discussions on $\dbH$}
\label{sect-LipH}

We first note that the $L$-integrability implies the Lipschitz property of $\dbH$ in $z$.

\begin{prop}
\label{prop-LipH}
Assume $\dbH$ is $L$-separable. Then $\dbH$ is Lipschitz continuous in $z$ with Lipschitz constant $L$.
\end{prop}
\proof Assume $\dbH = cl\{H_n, n\ge 1\}$ as in Definition \ref{defn-separable}. Fix $\th_i=(t,x, z_i)$, $i=1,2$. For any $y_1 \in \dbH(\th_1)$ and any $\e>0$, there exists $n_\e$ such that $|y_1 - H_{n_\e}(\th_1)|\le \e$. Then, by \reff{hreg},
\beaa
d\big(y_1, \dbH(\th_2)\big) \le d\big(H_{n_\e}(\th_1), \dbH(\th_2)\big) +\e \le \big|H_{n_\e}(\th_1)-H_{n_\e}(\th_2)\big| + \e\le L|z_1-z_2|+\e.
\eeaa
Since $\e>0$ is arbitrary, we have $d\big(y_1, \dbH(\th_2)\big) \le L|z_1-z_2|$ for all $y_1 \in \dbH(\th_1)$. Similarly $d\big(y_2, \dbH(\th_1)\big) \le L|z_1-z_2|$ for all $y_2 \in \dbH(\th_2)$. Then by \reff{dbHdistance} we obtain the desired Lipschitz property: $d\big(\dbH(\th_1), \dbH(\th_2)\big) \le L|z_1-z_2|$.
\qed

A natural and important question is: 
\bea
\label{Conjecture1}
\left.\ba{c}
\mbox{\it Assume $\dbH$ is uniformly Lipschitz continuous in $z$, then does \reff{VYH} hold true?}\\
\mbox{\it In particular,  is $\dbH$ separable or does Assumption \ref{assum-Lipschitz} hold true in this case?}
\ea\right.
\eea
We do not have a complete answer at this point. However, roughly speaking, if $\dbH$ can be parametrized with Lipschitz continuous parameters, then we have a positive answer. We illustrate this idea through the following example. 

\begin{eg}
\label{eg-parameter}
Let $B_{dN}$ denote the unit ball in $\dbR^{dN}$, and let $o: [0, T]\times \dbR^d \times \dbR^{dN} \to \dbR^N$ and $r: [0, T]\times \dbR^d \times \dbR^{dN} \times B_{dN} \to [0, \infty)$ be Lipschitz continuous in $z\in \dbR^{dN}$ with a Lipschitz constant $L_0$. Assume $r$  is continuous in $\zeta\in B_{dN}$. Then the following $\dbH$ is $2L_0$-separable:
\beaa
\dbH(t,x,z) := \big\{ o(t,x,z) + \iota r(t,x,z, \zeta)\zeta: \zeta\in B_{dN}, 0\le \iota\le 1 \big\}.
\eeaa 
\end{eg}
We note that $\dbH$ is a ball when $r$ is independent of $\zeta$. In general, it is connected but may not be convex. Moreover, by Lemma \ref{lem-union}, we may extend the result to unions of such sets and thus allow $\dbH$ to be not connected.

\ms
\proof Let $\{\zeta_n\}_{n\ge 1}$ be a dense subset of $B_{dN}$, and $\{\iota_m\}_{m\ge 1}$ be a dense subset of $[0, 1]$. Introduce the sequence of functions $H_{n,m}:\dbX\to \dbR^N$:
\beaa
H_{n,m}(\th) := o(\th) + \iota_m r(\th, \zeta_n) \zeta_n \in \dbH(\th),\q \th\in \dbX.
\eeaa
Then it is clear that $H_{n,m}$ is Lipschitz continuous in $z$ with Lipschitz constant $2L_0$. Moreover, for any $y:= o(\th) + \iota r(\th, \zeta)\zeta \in \dbH(\th)$ and any $\e>0$,  by the continuity of $r$ in $\zeta$, there exist $m, n$ such that 
\beaa
|r(\th, \zeta_n) - r(\th, \zeta)|\le \e,\q |\zeta_n - \zeta| \vee |\iota_m - \iota|\le {\e\over |r(\th, \zeta)|+1}.
\eeaa 
Then
\beaa
&&\dis \big|H_{n,m}(\th)-y\big| = \Big|\iota_m r(\th, \zeta_n) \zeta_n - \iota r(\th, \zeta)\zeta\Big|\\
&&\dis \le \big[|\iota_m - \iota| + |\zeta_n - \zeta|\big] r(\th, \zeta) +  |r(\th, \zeta_n) - r(\th, \zeta)| \le 3\e.
\eeaa
Since $\e$ is arbitrary, $\{H_{n,m}\}_{n,m\ge 1}$ is dense in $\dbH$, and hence $\dbH$ is $2L_0$-separable.
\qed

When $\dbH$ is a ball, a stronger version of Assumption \ref{assum-Lipschitz} can be verified straightforwardly.
\begin{eg}
\label{eg-parameter2}
Consider the setting in Example \ref{eg-parameter} with $r = r(t,x,z)$ independent of $\zeta$. Let $\eta\in \dbL^2(\dbF; \dbR^N)$ and $Z\in \dbL^2(\dbF; \dbR^{dN})$. Since $\dbF=\dbF^X$, by abusing the notation we write $\eta_t = \eta(t, X_\cd)$ for a deterministic mapping $\eta: [0, T]\times C([0, T]; \dbR^d)\to \dbR^N$. Then the following function $H$, which projects $\eta$ onto $\dbH$, satisfies Assumption \ref{assum-Lipschitz} with $L=4L_0$ and $\e=0$:
\beaa
H(t, \bx, z) := \left\{\ba{lll} 
\dis \eta(t, \bx),\q \mbox{if}~ |\eta(t, \bx) - o(t,\bx_t,z)| \le r(t,\bx_t,z);\\ 
\dis o(t,x,z) + {r(t,\bx_t,z)\over |\eta(t, \bx) - o(t,\bx_t,z)|} \big[\eta(t, \bx) - o(t,\bx_t,z)\big],\q \mbox{otherwise}.
\ea\right.
\eeaa
\end{eg}
\proof First, by the definition of $\dbH$ and $H$  it is clear that $H(t,\bx, z)\in \dbH(t,\bx_t,z)$ and $\big|\eta(t,\bx) -  H(t, \bx, z)\big| = d\big( \eta(t,\bx),  \dbH(t, \bx_t, z)\big)$. Then it suffices to verify that $H$ is $4L_0$-Lipschitz continuous in $z$.  We prove it in three cases. Fix $z_1, z_2$.

{\it Case 1}. $\eta(t,\bx)$ is in both $\dbH(t,\bx_t, z_i)$. Then clearly $|H(t,\bx, z_1) - H(t,\bx,z_2)| =0$. 

{\it Case 2}. $\eta(t,\bx)$ is in neither $\dbH(t,\bx_t, z_i)$. Then by the Lipschitz property of $o$ and $r$ in $z$ and noting that $|\eta(t, \bx) - o(t,\bx_t,z_i)| > r(t,\bx_t,z_i)$, it is straightward to see that $ |H(t,\bx, z_1) - H(t,\bx,z_2)| \le 4L_0|z_1-z_2|$.

{\it Case 3}. $\eta(t,\bx)$ is in exactly one of $\dbH(t,\bx_t, z_i)$, and assume without loss of generality that it's in  $\dbH(t,\bx_t, z_1)$. Then, denoting $\hat \eta(t,\bx, z):= \eta(t,\bx) - o(t,\bx_t,z)$, 
\beaa
&&\dis \big|H(t,\bx, z_1) - H(t,\bx, z_2)\big| = \Big| \hat\eta(t,\bx,z_2) - {r(t,\bx_t,z_2)\over |\hat\eta(t, \bx,z_2)|} \hat\eta(t, \bx,z_2)\Big|\\
&&\dis =\big|\hat\eta(t, \bx,z_2)\big| - r(t,\bx_t,z_2)   \le  \big|\hat\eta(t, \bx,z_1)\big| + L_0|z_1-z_2| - r(t,\bx_t,z_2)\\
&&\dis \le r(t,\bx_t,z_1) - r(t,\bx_t,z_2)  + L_0|z_1-z_2| \le 2L_0|z_1-z_2|.
\eeaa
Put together, $H$ is $4L_0$-Lipschitz continuous in $z$.
\qed

\section{Appendix}
\label{sect-appendix}
\no{\bf Proof of Proposition \ref{prop-HH0}.}
 It is clear that $\dbH_0 \subset \dbH$. To see the opposite inclusion, we fix $\th\in \dbX$ and  $y\in\dbH(\th)$. For any $\e>0$, by \reff{olYa*} there exists $a^\e\in \cE_\e(\th)$ such that 
   \bea
  \label{ancompact}
  |h(\th,a^\e)-y|\leq\e,\q\mbox{and}\q
h_i(\th,a^{\e})\geq \ol h_i(\th, a^{\e})-\e_n, \q  i=1,\cds, N.
  \eea
  Since $A$ is compact, there exists a sequence $\e_n\downarrow 0$ such that $a^{\e_n}\to a^*\in A$. Then, by applying $\e_n$ in \reff{ancompact} and sending $n\to\infty$, it follows from \reff{hreg}  that
 \beaa
 y= h(\th, a^*),\q\mbox{and}\q h_i(\th,a^*)\geq \ol h_i(\th, a^{*}), \q i=1,\cds, N.
  \eeaa
  That is, $a^*\in \cE(\th)$ and $y =  h(\th,a^*) \in \dbH_0(\th)$. Since $y\in \dbH(\th)$ is arbitrary, we obtain $\dbH(\th)\subset \dbH_0(\th)$, and hence equality holds. 
\qed

\bs
\no{\bf Proof of Proposition \ref{prop-VH0}.} We first prove the ``if'' direction. Assume \reff{etaH0} holds true and denote $\Th_t := (t, X_t, Z^\eta_t)$. Since by \reff{hreg} $h$ and $\ol h$ are continuous in $a$, one can easily construct a measurable process $\a^*$ such that  $\eta^i_t = h_i(\Th^i_t, \a^*_t)= \ol h_i(\Th^i_t, \a^*_t)$, $dt\times d\dbP$-a.s. Then $Y^{\a^*,i}_0= Y^{\eta,i}_0= \ol Y^{\a^*,i}_0$. By \reff{olYa*},
this implies that $\a^*\in \cE$, and thus $y = Y^{\eta}_0 = J(\a^*) \in \dbV_0$.

We next prove the ``only if'' direction. Assume $y = Y^{\a^*}_0\in \dbV_0$ for some $\a^*\in \cE$. Denote $\Th_t := (t, X_t, Z^{\a^*}_t)$, $\ol\Th_t := (t, X_t, \ol Z^{\a^*}_t)$, $\D Y_t := \ol Y^{\a^*}_t - Y^{\a^*}_t$, $\D Z_t := \ol Z^{\a^*}_t - Z^{\a^*}_t$. Then we have
\beaa
\D Y^i_t = \int_t^T \Big[\ol h_i( \ol \Th^i_s, \a^{*}_s) - h_i(\ol \Th^i_s, \a^*_s) + \b^i_s \D Z^i_s\Big] ds- \int_t^T \D Z^i_s \cd dB_s,
\eeaa
for some appropriate bounded process $\b^i$. By \reff{olYa*} and recalling \reff{Mbeta}, we have, 
\beaa
0 = \D Y^i_0 = \dbE\Big[M^{\b^i}_T \int_0^T \big[ ~\!\ol h_i(\ol\Th^i_t, \a^{*}_t) - h_i(\ol\Th^i_t, \a^*_t) \big]dt\Big].
\eeaa
Note that $\ol h_i(\ol\Th^i_t, \a^{*}_t) - h_i(\ol\Th^i_t, \a^*_t)\ge 0$, $dt\times d\dbP$-a.s. Then we must have equality:  $\ol h(\ol\Th_t, \a^{*}_t) = h(\ol\Th_t, \a^*_t)$. This implies $\a^*_t \in \cE_0(\ol\Th_t)$ and $\eta_t := h(\ol\Th_t, \a^*_t) \in \dbH_0(\ol\Th_t)$, $dt\times d\dbP$-a.s. It is clear that $(\ol Y^{\a^*}, \ol Z^{\a^*}) = (Y^\eta, Z^\eta)$, we thus obtain \reff{etaH0}.
\qed

\bs
\no{\bf Proof of Proposition \ref{prop-construction}.} We proceed in two steps. Fix $\eta, Z, \e, \Th$ as required. 

{\bf Step 1.} We first assume $\eta$ and $Z$ are continuous in $t$, and $|Z|\le R$ for some constant $R>0$, a.s. By \reff{hreg} $h$ and hence $\ol h$ are uniformly continuous in all variables $(\th, a)$ when $|z|\le R$, with a modulus of continuity function $\rho_R$. Let $\d>0$ be a small number which will be specified later.  Define $\t_0 :=0$, and for $n=0,1,\cds$,
\bea
\label{taun}
\t_{n+1} := \Big\{t\ge \t_n:  |\Th_t - \Th_{\t_n}| + d\big(\dbH(\Th_t), \dbH(\Th_{\t_n})\big)+|\eta_t - \eta_{\t_n}| \ge \d\Big\}\wedge T.
\eea 
By the continuity of $X, Z,\eta$, we have $\t_n = T$ for $n$ large enough, a.s.  Next,  let  $\{D_m\}_{m\ge 1}$ be a partition of $\{(\th, a)\in \dbX\times A: |z|\le R\}$ such that, for any $m\ge 1$ and for arbitrarily fixed $(\th_m, y_m)\in cl(D_m)$,
\bea
\label{partition}
|\th-\th_m|+  d(\dbH(\th), \dbH(\th_m)) + |y-y_m| \le \d,\q \forall (\th, y)\in D_m.
\eea
Since $\dbH(\th_m)$ is compact, there exists $y_m^*\in \dbH(\th_m)$ such that $|y_m-y_m^*| = d\big(y_m, \dbH(\th_m)\big)$. Moreover, since $y_m^*\in \dbH(\th_m)\subset \dbH_\d(\th_m)$, there exists $a_m\in \cE_\d(\th_m)$ such that $|y_m^* - h(\th_m, a_m)|< \d$.  We now construct $\a$ as follows  which is obviously in $\cA$:
\bea
\label{constructa}
\a_t = \sum_{n=0}^\infty \1_{[\t_n, \t_{n+1})}(t) \sum_{m=1}^\infty a_m \1_{D_m}(\Th_{\t_n}, \eta_{\t_n}).
\eea
We shall verify that $\a$ satisfies all the requirements when $\d>0$ is small enough. 

Fix $(t, \o)$ and the corresponding $n,m$ such that $t\in [\t_n, \t_{n+1})$ and $(\Th_{\t_n}, \eta_{\t_n})\in D_m$. First, since $a_m\in \cE_\d(\th_m)$, by \reff{olYa*} we have $h_i(\th^i_m, a_m) \ge \ol h_i(\th^i_m, a_m) - \d$ for all $i=1,\cds, N$. Then, at $(t,\o)$, by \reff{partition} we have
\beaa
 \dis 0 &\le& \ol h_i(\Th_t^i, \a_t) - h_i(\Th_t^i, \a_t) = \ol h_i(\Th_t^i, a_m) - h_i(\Th_t^i, a_m) \le 2 \rho_R\big(|\Th_t - \th_m|\big) + \d\\
 &\le& 2 \rho_R\Big(\big|\Th_t- \Th_{\t_n}\big| + \big|\Th_{\t_n}-\th_m\big|\Big) + \d \le 2 \rho_R(2\d)+\d.
 \eeaa
 Thus, $\a_t \in \cE_\e(\Th_t)$ when $2 \rho_R(2\d)+\d< \e$.
 
 Next, for any $y_0\in \dbR^N$ and $\th, \tilde \th \in \dbX$, note that
 \beaa
&&\dis d\big(y_0, \dbH(\th)\big) - d\big(y_0, \dbH(\tilde\th)\big)  = \inf_{y\in \dbH(\th)} |y_0-y| - \inf_{\tilde y\in \dbH(\tilde \th)} |y_0-\tilde y|\\
&&\dis = \sup_{\tilde y\in \dbH(\tilde \th)}  \inf_{y\in \dbH(\th)} \big[|y_0-y| - |y_0-\tilde y|\big] \le  \sup_{\tilde y\in \dbH(\tilde \th)}  \inf_{y\in \dbH(\th)} |y-\tilde y| =  \sup_{\tilde y\in \dbH(\tilde \th)}  d\big(\tilde y, \dbH(\th)\big).
 \eeaa
 Similarly, we have $d\big(y_0, \dbH(\tilde\th)\big) - d\big(y_0, \dbH(\th)\big)  \le \sup_{y\in \dbH(\th)}  d\big(y, \dbH(\tilde\th)\big)$.
Then 
\bea
\label{dHreg}
\big|d\big(y_0, \dbH(\th)\big) - d\big(y_0, \dbH(\tilde\th)\big)\big| \le d\big(\dbH(\th), \dbH(\tilde\th)\big).
\eea
Thus
 \beaa
&&\dis \Big|\big|\eta_t - h(\Th_t, \a_t)\big| - d\big(\eta_t, \dbH(\Th_t)\big)\Big|\\
&&\dis \le \Big|\big|\eta_t - h(\Th_t, \a_t)\big| - |y_m-y_m^*|\Big| + \Big| d\big(y_m, \dbH(\th_m)\big)-d\big(\eta_t, \dbH(\Th_t)\big)\Big|\\
&&\dis \le 2|\eta_t  - y_m| + |y_m^* -  h(\Th_t, \a_t)| +  \Big| d\big(\eta_t, \dbH(\th_m)\big)-d\big(\eta_t, \dbH(\Th_t)\big)\Big| \\
&&\dis \le 2|\eta_t  - y_m| + |y_m^* -  h(\Th_t, \a_t)| + d\big(\dbH(\th_m), \dbH(\Th_t)\big) \\
&&\dis \le 2\big[|\eta_t - \eta_{\t_n}| + |\eta_{\t_n} - y_m|\big] + |y_m^* - h(\th_m, a_m)| + | h(\th_m, a_m) - h(\Th_{\t_n}, a_m)\big| \\
&&\dis\q + | h(\Th_{\t_n}, a_m)- h(\Th_t, a_m)\big| + d\big(\dbH(\th_m), \dbH(\Th_{\t_n})\big)+d\big(\dbH(\th_{\t_n}), \dbH(\Th_t)\big)\\
&&\dis \le C\d + 2\rho_R(\d).
 \eeaa
 Then, by setting $Z^\e:=Z$,  clearly $(\a, Z)$ satisfy all the requirements when $\d$ is small enough.
 
 {\bf Step 2.} We next consider general $\eta$ and $Z$. Fix $R>0$ large and $\d>0$ small. Let $\eta^R := \eta \1_{\{|\eta|\le R\}}$ and $Z^R:=  Z \1_{\{|Z|\le R\}}$ denote the truncations. By standard arguments there exist continuous processes $\eta^{R,\d}$, $Z^{R,\d}$ such that they are bounded by $R$ and  
 \bea
 \label{ZRd}
 \dbE\Big[\int_0^T \big[|\eta^{R,\d}_t-\eta^R_t|^2 + |Z^{R,\d}_t-Z^R_t|^2\big]dt\Big]\le \d.
 \eea
 Denote  $\Th^{R,\d}_t := (t, X_t, Z^{R,\d}_t)$. By Step 1, there exists $\a$ such that  
 \bea
 \label{etaRd}
\a_t \in \cE_\d(\Th^{R,\d}_t),\q \Big|\big|\eta^{R,\d}_t - h(\Th^{R,\d}_t, \a_t)\big| - d\big(\eta^{R,\d}_t, \dbH(\Th^{R,\d}_t)\big)\Big| \le \d.
\eea
We now verify that $(\a, Z^{R,\d})$ satisfy all the requirements when $R$ is large and $\d$ is small.

First, by \reff{ZRd} we have
\bea
\label{Zest}
\dis \dbE\Big[\int_0^T|Z^{R,\d}_t-Z_t|^2dt\Big] &\le& C\dbE\Big[\int_0^T |Z^R_t-Z_t|^2dt\Big]+C\d \nonumber\\
\dis &=& C\dbE\Big[\int_0^T |Z_t|^2\1_{\{|Z_t|\ge R\}}dt\Big]+C\d.
\eea
Clearly the right side will be less than $\e$ when $\d\le {\e\over 2C}$ and $R$ is sufficiently large, which may depend on $Z$. This, together with \reff{etaRd}, verifies the first line of \reff{construction}.

To see the second line to \reff{construction},  by  \reff{etaRd}, \reff{hreg}, and \reff{dHreg},  we have
\beaa
\D &:=&\dbE\Big[\int_0^T \Big|\big|\eta_t - h(\Th_t, \a_t)\big| - d\big(\eta_t, \dbH(\Th_t)\big)\Big|^2dt\Big] \\
 &\le& C \dbE\Big[\int_0^T \Big[\big|\eta_t - \eta^{R,\d}_t\big|^2 +  \big|Z^{R,\d}_t-Z_t\big|^2 + d^2\big(\dbH(\Th^{R,\d}_t), \dbH(\Th_t)\big)\Big]dt\Big]+C\d\\
 &\le& C \dbE\Big[\int_0^T \Big[|\eta_t|^2\1_{\{|\eta_t|\ge R\}} + |Z_t|^2\1_{\{|Z_t|\ge R\}}  + d^2\big(\dbH(\Th^{R,\d}_t), \dbH(\Th_t)\big)\Big]dt\Big]+C\d,
\eeaa
where the last inequality is due to \reff{Zest} and a similar estimate for $\eta$. Next, denote $X^R_t:= X_t\1_{\{|X_t|\le R^2\}}$, $\Th^R_t:=(t, X^R_t, Z^R_t)$, and $\tilde \Th^{R,\d}_t:=(t, X^R_t, Z^{R,\d}_t)$. Then, by Proposition \ref{prop-VHcompact} (iii),
\beaa
&&\dis  d\big(\dbH(\Th^{R,\d}_t), \dbH(\Th_t)\big)  \le  d\big(\dbH(\Th^{R,\d}_t), \dbH(\tilde \Th^{R,\d}_t)\big) + d\big(\dbH(\tilde\Th^{R,\d}_t), \dbH(\Th^{R}_t)\big) + d\big(\dbH(\Th^{R}_t), \dbH( \Th_t)\big)\nonumber\\
&&\dis\qq \le CR\1_{\{|X_t|\ge R^2\}} + d\big(\dbH(\tilde\Th^{R,\d}_t), \dbH(\Th^{R}_t)\big) + C[1+|Z_t|] [\1_{\{|X_t|\ge R^2\}} +\1_{\{ |Z_t|\ge R\}}\big].
\eeaa
Note that $R\1_{\{|X_t|\ge R^2\}}\le {|X_t|\over R}$. Then we have
\beaa
 &&\dis \D \le \dbE\Big[\int_0^T d^2\big(\dbH(\tilde\Th^{R,\d}_t), \dbH(\Th^{R}_t)\big) dt\Big] + C\d\\
 &&\dis\q + C \dbE\Big[\int_0^T \Big[|\eta_t|^2\1_{\{|\eta_t|\ge R\}} + [1+|Z_t|^2]  [\1_{\{|X_t|\ge R^2\}} +\1_{\{ |Z_t|\ge R\}}\big] +  {1\over R^2}|X_t|^2\Big]dt\Big].
\eeaa
Choose $R$ large enough such that the second line above is less than ${\e\over 2}$. Then
\beaa
 \D \le \dbE\Big[\int_0^T d^2\big(\dbH(\tilde\Th^{R,\d}_t), \dbH(\Th^{R}_t)\big) dt\Big] + C\d+{\e\over 2}.
\eeaa
 Moreover, since $\dbH$ is continuous and hence uniformly continuous on $\{\th\in\dbX: |x|\le R^2, |z|\le R\}$, then for fixed $R$, by \reff{ZRd} we have $\dis\lim_{\d\to 0}\dbE\Big[\int_0^T d^2\big(\tilde\Th^{R,\d}_t), \dbH(\Th^{R}_t)\big)dt\Big]=0$. Thus $\D \le \e$ when $\d>0$ is small enough. This proves the second line of \reff{construction}.
\qed

\bibliographystyle{siam}
\bibliography{refs}

\end{document}